\documentclass[12pt,a4paper,titlepage]{article}
\usepackage{mathtext,amssymb,amsmath,latexsym,amsthm}
\usepackage[T1,T2A]{fontenc}
\usepackage[english]{babel}
\usepackage[dvips]{graphicx}

\font\got = eufm10  scaled \magstep2
 
\begin{document}

\begin{Large}
 \centerline{\bf A spectrum associated with}
\centerline{\bf  Minkowski diagonal continued fraction}
\end{Large}
\vskip+0.5cm
\centerline{ by {\bf
Alena Aleksenko}\footnote{Department of Mathematics, Aveiro
University, Aveiro 3810, Portugal.   This work was supported by {\it
FEDER} founds through {\it COMPETE}--Operational Programme Factors of
Competitiveness (``Programa Operacional Factores de Competitividade'')
and by Portuguese founds through the {\it Center for Research and
Development in Mathematics and Applications} (University of Aveiro) and
the Portuguese Foundation for Science and Technology
(``FCT--Fund\c{c}\~{a}o para a Ci\^{e}ncia e a Tecnologia''), within
project PEst-C/MAT/UI4106/2011 with COMPETE number
FCOMP-01-0124-FEDER-022690}}
\vskip+1.0cm

Let  $\alpha $
be real irrational number.
The function 
$\mu_\alpha (t)$
is defined as follows.

The Legendre theorem states that 
 if 
\begin{equation}\label{1}
\left|
\alpha - \frac{A}{Q}\right|< \frac{1}{2Q^2},\,\,\,
(A,Q) =1
\end{equation}
then the fraction $\frac{A}{Q}$  is a convergent fraction for the continued fraction expansion of
$\alpha$. The converse statement is not true.
It may happen that
$\frac{A}{Q}$ is a convergent to $\alpha$
but (\ref{1}) is not valid.
One should consider the sequence of the denominators of the convergents to $\alpha$
for which (\ref{1}) is true.
Let this sequence be
$$
Q_0<Q_1<\cdots < Q_n<Q_{n+1}<\cdots .
$$
Then 
for $\alpha \not \in \mathbb{Q}$ the function
$\mu_\alpha (t) $ is defined by
$$
\mu_\alpha (t) =
 \frac{Q_{n+1}-t}{Q_{n+1}-Q_n}\cdot ||Q_n \alpha||
+
 \frac{t-Q_n}{Q_{n+1}-Q_n}
\cdot ||Q_{n+1} \alpha||
,\,\,\
Q_n \le t \le Q_{n+1}.
$$
From the other hand,
for every $\nu$ one of the consecutive convergent fractions $\frac{p_\nu}{q_\nu},\,\frac{p_{\nu+1}}{q_{\nu+1}}$
to $\alpha$  satisfies (\ref{1}). 
So either
$$
(Q_n, Q_{n+1}) = (q_\nu, q_{\nu+1})
$$
for some $\nu$,
 or
$$
(Q_n, Q_{n+1}) = (q_{\nu-1}, q_{\nu+1})
$$
for some $\nu$.

Actually the function $\mu_\alpha (t)$ was considered by Minkowski
\cite{MI}.  There exists an alternative geometric definition of $\mu_\alpha (t)$.
Some related facts were discussed in \cite{1,2}.

The quantity
$$
\hbox{\got m} (\alpha) = \limsup_{t\to +\infty} t\cdot \mu_\alpha (t).
$$
was considered in \cite{2}. 
An explicit formula for the value of $
\hbox{\got m} (\alpha)$ in terms of continued fraction expansion for $\alpha$
was proved in \cite{2}. It is as follows.
Put
\begin{equation}\label{em}
\hbox{\got m}_n (\alpha)=
\begin{cases}
 G(\alpha_\nu^*,\alpha_{\nu+2}^{-1}),\,\,\,\text{if}\,\,\, (Q_n,Q_{n+1}) =(q_{\nu-1},q_{\nu+1})\,\,\,\text{with some } \,\,\,\nu,
\cr
F(\alpha_{\nu+1}^*,
\alpha_{\nu+2}^{-1}),\,\,\,\text{if}\,\,\, (Q_n,Q_{n+1}) =(q_{\nu},q_{\nu+1}) \,\,\,\text{with some } \,\,\,\nu,
\end{cases}
 \end{equation}
where
$$
G(x,y) = \frac{x+y+1}{4}
,\,\,\,\,\,
F(x,y) = \frac{(1 - xy)^2}{4(1+xy) (1-x)(1-y)} 
$$ and
$\alpha_\nu, \alpha_\nu^*$ 
come from continued fraction expansion to 
$$
\alpha = [a_0;a_1,a_2,\dots,a_t,\dots]
$$
in such a way:
$$
\alpha_\nu = [a_\nu;a_{\nu+1},...],\,\,\,\, \alpha_\nu^* = [0;a_\nu, a_{\nu-1},...,a_1].
 $$
Then
$$
\hbox{\got m} (\alpha) = \limsup_{n\to +\infty} \hbox{\got m}_n (\alpha),
$$
The specrtum
 $$
\mathbb{M} = \{ m \in \mathbb{R}: \,\,\, \exists \alpha \in \mathbb{R}\setminus\mathbb{Q} \,\,
\text{such that}\,\, m = \hbox{\got m} (\alpha)\}.
  $$
was studied in \cite{2}.
It was proven there that $\mathbb{M}\subset \left[\frac{1}{4},\frac{1}{2}\right]$
and that $\frac{1}{4},\frac{1}{2} \in \mathbb{M}$.
However no further structure of the spectrum $\mathbb{M}$ is known.

In this paper we consider the spectrum
$$
\mathbb{I} =
\{ m\in \mathbb{R}:\,\,\,\exists \alpha \in \mathbb{R}\setminus\mathbb{Q} \,\,\text{such that}\,\,
\hbox{\got i} (\alpha) = m\},
$$
where
$$
\hbox{\got i} (\alpha ) =\liminf_{n \to \infty} \hbox{\got m}_n (\alpha )
$$
(however, compared to $\hbox{\got m}(\alpha)$, this quantity has no clear Diophantine sense).
  
It is clear that
$$
\min \mathbb{I} = \frac{1}{4},\,\,\,\,\,
\max \mathbb{I} = \frac{1}{2}.
$$

{\bf Theorem.} \,\,{\it There exists positive $\omega_0$ such that
 $$
 \left[\frac{1}{4}, \omega_0\right] \subset \mathbb{I}
.
$$
}

The proof is based ol M. Hall's ideas (see \cite{hall}).
It uses technique from \cite{MoHH}.
 
{\bf Remark.}\,\,{\it
An explicit formula for $\omega_0$
 may be obtained from the proof below.
It is interesting to get optimal estimates for the value of $\omega_0$.}
 
We need some well known results.

Recall 
the definition of a $\tau$-set $ {\cal F} \subset \mathbb{R}$.
The set ${\cal F}$ must be of the form
$$
{\cal  F} = {\cal S}\setminus \left(\bigcup_{\nu+1}^\infty \Delta_\nu\right),
$$
where  ${\cal S}\subset \mathbb{R}$ is  a segment, and $ \Delta_\nu\subset{\cal S},\,\, \nu  =1,2,3,...$ is an ordered sequence of 
disjoint intervals.
Moreover for every $t$ 
if
$$
{\cal S }\setminus \left(\bigcup_{\nu+1}^{t-1} \Delta_\nu \right)
= \bigcup_{j=1}^r {\cal M}_j
$$
is a union of segments ${\cal M}_j$
and
$\Delta_t \subset {\cal M}_{j^*}
$
then
$$
{\cal M}_{j^*} = {\cal N}^1 \sqcup \Delta_t \sqcup {\cal N}^2,
$$
and 
$$
\min (|{\cal N}^1|,  |{\cal N}^2|
 )
\ge
\tau |\Delta_t|.
$$

Consider the set 
 $$
{\cal 
F}_5 =
\{
\alpha = [0; b_1, b_2, b_3, ... ]:\hspace{2mm}
b_\nu  \le 5  \,\,\,\,\forall\,\, \nu
\}
$$
consisting of all irrational real numbers from the unit interval $(0,1)$ with partial quotients bounded by $5$.
One can easily see that 
\begin{equation}\label{abe}
\begin{cases}
A =
\min
{\cal F}_5
=
[0;\overline{5,1}]
=
\frac{\sqrt{45}-5}{10}
=0.1708^+,\cr
B = \max
{\cal F}_5
=
[0;\overline{1,5}]
=
\frac{\sqrt{45}-5}{2} = 0.85410^+.
\end{cases}
\end{equation}
Put
$$
{\cal S}_5 =
[A,B] \subset [0,1].
$$

The following lemma comes from the results of the papers \cite{asti} or \cite{psia}. 

{\bf Lemma 1.}\,\,{\it The set ${\cal F}_5$ is a $\tau$-set with
$ \tau = \tau_5 = 1.788^+$.}

Let
$ H(x,y) : \hspace{2mm}
{\cal S} \times {\cal S} \to \mathbb{R}$
be a function in two variables of the class $ G \in C^1 ({\cal S}\times {\cal S})$.
Consider the set
$$ {\cal J} = \{ z \in \mathbb{R} :\hspace{2mm}
\exists\, x,y \in {\cal S}\hspace{2mm}  z = H(x,y) \}.
$$
By continuousity argument  ${\cal J}$ is a segment.

{\bf  Lemma 2.}\,\,{\it
Suppose that the derivatives 
$\partial H/\partial x, \partial H/\partial y$
do not take zero values on the box $ {\cal S}\times {\cal S}$.
 Suppose that ${\cal F}$ is $\tau$-set and  ${\cal S} = [\min {\cal F}, \max {\cal F}]$.
Suppose that
\begin{equation}\label{ouslovie}
\tau 
\ge  \max_{x,y \in {\cal S}}
\,\,
\max
\left(\left|
\frac{\partial H/\partial x}{\partial H/\partial y}\right|,
\left|
\frac{\partial H/\partial y}{\partial H/\partial x}\right|\right)
.
\end{equation}
 Then
$$
\{ z:\,\,\,
\exists\, x,  y \in {\cal F}\,\,\,\text{ such that }\,\,\, z = H(x,y
) \}= {\cal J}.
$$
}

Lemma 2  is a staightforward 
generalization of a result from \cite{MoHH}.
 We do not give its proof here as the proof follows the argument from \cite{MoHH} word-by-word.

Now we are able to conclude the proof of Theorem.

We consider pairs of integers $(R_1,R_2)$ of  the form
\begin{equation}\label{reeva}
(R_1,R_2)\,\, =
\,\,
(R,R)\,\,\,\,
\text{or}\,\,\,\,
(R,R+1)
\end{equation}
with $R\ge 6$.
Consider a function
$$
H_{R_1,R_2} (x,y) =
F\left(
\frac{1}{R_1+x}, \frac{1}{R_2+y}\right).
$$
For $R_1,R_2$ under consideration the function $
H_{R_1,R_2} (x,y)$ decreases both in $x$ and in $y$.

For $0<x,y<1$ put
$$
\varphi (x,y) =
(1-3x+3xy - x^2y)(1-x).
$$
For any $y \in (0,1) $ the function 
$\varphi (x,y)$  decreases in $x$.
For any $ x \in (0,1) $ the function 
$\varphi (x,y)$  increases in $y$.
Now
$$
\frac{\partial F/\partial y}{\partial F/\partial x}
=\frac{\varphi (x,y)}{\varphi (y,x)},
$$
and
$$
\left|
\frac{\partial H_{R_1,R_2}/\partial y}{\partial H_{R_1,R_2}/\partial x}\right|
=\frac{ \varphi\left(\frac{1}{R_1+x}, \frac{1}{R_2+y}\right)}
{ \varphi\left(\frac{1}{R_2+y}, \frac{1}{R_1+x}\right)}
\,\left(
\frac{R_1+x}{R_2+y}
\right)^2.
$$
Easy calculation shows that  for $R_1, R_2 \ge 6$
one has
$$
\max_{x,y \in {\cal S}_5}
\,\,
\max
\left(\left|
\frac{\partial H_{R_1,R_2}/\partial x}{\partial H_{R_1,R_2}/\partial y}\right|,
\left|
\frac{\partial H_{R_1,R_2}/\partial y}{\partial H_{R_1,R_2}/\partial x}\right|\right)
= 
 $$
$$
=
\frac{ \varphi\left(\frac{1}{R_1+B}, \frac{1}{R_2+A}\right)}
{ \varphi\left(\frac{1}{R_2+A}, \frac{1}{R_1+B}\right)}
\,\left(
\frac{R_1+B}{R_2+A}
\right)^2
\le
\frac{ \varphi\left(\frac{1}{R_1+B}, \frac{1}{R_1+A}\right)}
{ \varphi\left(\frac{1}{R_1+A}, \frac{1}{R_1+B}\right)}
\,\left(
\frac{R_1+B}{R_1+A}
\right)^2
\le
$$
$$
\le
\frac{ \varphi\left(\frac{1}{6+B}, \frac{1}{6+A}\right)}
{ \varphi\left(\frac{1}{6+A}, \frac{1}{6+B}\right)}
\,\left(
\frac{6+B}{6+A}
\right)^2
 = 1.363^+< \tau_5.
$$
Here $A$ and $B$ are defined in (\ref{abe})
and in the last inequalities we use the bounds $6 \le  R_1\le R_2$ 
 which follows from (\ref{reeva}).

We see that for any $R_1, R_2$ under consideration and for $\tau_5$-set ${\cal F}_5$ the condition
(\ref{ouslovie}) is satisfied. We apply Lemma 2 to see that the
image of the set ${\cal F}_5\times {\cal F}_5$ under the mapping
$H_{R_1,R_2} (x,y)$ is just the segment
$$
{\cal J}_{R_1,R_2} =
[ H_{R_1,R_2} (B,B), H_{R_1,R_2}(A,A)].
$$
But
$$
H_{R,R}(B,B) < H_{R,R+1} (A,A)
$$ and
$$
H_{R,R+1}(B,B) < H_{R+1,R+1} (A,A).
$$
That is why
if we put
$$
\omega_0 = H_{R_0,R_0} (A,A).
$$
with $R_0 \ge 6$ we get
$$
\bigcup_{R\ge R_0} {\cal J}_{R,R}\,
\cup\,\bigcup_{R\ge R_0} {\cal J}_{R, R+1} =
(1/4, \omega_0]
.
$$

Take $ m \in (0,\omega_0]$. Then there exists $R_1,R_2$ such that
$$
m \in 
{\cal J}_{R_1,R_2}$$
and
there exist 
$$
 \beta = [0;b_1,b_2,...,b_\nu,...],\,\,\,
\gamma = [0; c_1,c_2,...,c_\nu,...],\,\,\,
 \beta, \gamma \in {\cal F}_5,
$$
such that
$$
F\left( \frac{1}{R_1+\alpha}, \frac{1}{R_2+\beta}\right) = m.
$$
Now we take
$$
\alpha =
[0;\underbrace{a_1,R_1,R_2,b_1}_1,\underbrace{a_2,a_1,R_1,R_2,b_1,b_2}_2,...,
\underbrace{a_\nu, a_{\nu-1},...,a_2,a_1, R_1,R_2,b_1,b_2,...,b_{\nu-1},b_\nu}_\nu,...
].
$$
Standard argument shows that 
for $n_\nu$ defined from
$$
\frac{p_{n_\nu}}{q_{n_\nu}}=
[0;{a_1,R_1,R_2,b_1},{a_2,a_1,R_1,R_2,b_1,b_2},...,
a_\nu, a_{\nu-1},...,a_2,a_1, R_1]
$$
one has
$$
\lim_{\nu \to +\infty} F(\alpha_{n_\nu}^*, \alpha_{n_\nu+1}^{-1}) = m
.$$
At the same time for $F(\cdot, \cdot)$ and $G(\cdot, \cdot )$  we have
$$
\inf_{n \neq n_\nu \,\forall \nu}  F(\alpha_n^*,\alpha_{n+1} )   >\omega_0 
$$
and
$$
\inf_{n \in \mathbb{Z}_+} G(\alpha_n^*,\alpha_{n+2}^{-1} )  >\omega_0 
,$$
for large $R_0$.
So 
$\hbox{\got i}(\alpha ) = m$ and everything is proved.$\Box$

\end{document}